\documentclass[a4paper,11pt]{article}
\usepackage[english]{babel}
\usepackage{graphicx}
\usepackage{bbding}
\usepackage[latin1]{inputenc}
\usepackage{fancyhdr}
\usepackage{amsmath,amssymb,color,epsfig}
\usepackage{mathrsfs}
\usepackage{eufrak}
\usepackage{dsfont}
\usepackage{upgreek}
\usepackage{fancyhdr}
\usepackage[dvipsnames]{xcolor}
\newtheorem{theorem}{Theorem}[section]
\newtheorem{proposition}[theorem]{Proposition}
\newtheorem{corollary}[theorem]{Corollary}
\newtheorem{remark}[theorem]{Remark}
\newtheorem{lemma}[theorem]{Lemma}

\def\1{\mathds{1}}

\title{On some complete monotonic functions}
\author{Mohamed Bouali }                 %<-------------------
\date{}
\begin{document}

\maketitle
\begin{abstract} Motivated by open questions in the papers " Refinements and sharpenings of some double inequalities for
bounding the gamma function" and "Complete monotonicity and monotonicity of two functions defined by two derivatives of a function involving trigamma function",we confirm among other results and disprove other one.
\end{abstract}
\section{Introduction}
Completely monotonic functions have attracted the attention of many authors. Mathematicians
have proved many interesting results on this topic. For example, Koumandos \cite{kou1} obtained upper and
lower polynomial bounds for the function $x/(e^x -1)$, $x > 0$, with coefficients of the Bernoulli numbers
$B_k$. This enabled him to give simpler proofs of some results of H. Alzer and F. Qi et al., concerning
complete monotonicity of certain functions involving the functions $\Gamma(x)$, $\psi(x)$ and the polygamma
functions $\psi^{(n)}(x)$, $n = 1, 2, ....$, \cite{clar}.

A function $f$ is said to be completely monotonic on an interval $I$ if $f$ has derivatives of all
orders on $I$ which alternate successively in sign, that is, $(-1)^nf^{(n)}(x)\geq 0$ for all $x\in I$ and all $n\in\Bbb N\cup\{0\}$. %The function is said to be strictly completely monotonic, if the inequalities on the derivatives are all strict.
See for example [\cite{mit}, Chap VIII], [\cite{sch}, Chap I], and [\cite{wid}, Chap IV].

A notion of logarithmically completely monotonic is introduced in reference \cite{Qi10} and \cite{Qi11}. A positive function $f$ is said to be logarithmically completely monotonic on an open
interval $I$, if $f$ satisfies $(-1)^n(\log f(x))^{(n)}\geq 0$ for all $x\in I$ and all $n\in\Bbb N$.
 %When the inequalities are strict we say that $f$ is strictly Logarithmically completely monotonic.
For more informations on completely monotonic functions, see for instance \cite{guo1,kou2,kou3,kou4,kou5,kou6, Qi1, Qi2, Qi3, Qi4} and closely related references.

%Let $f(x)$ be a completely monotonic function on $(0,+\infty)$ and denote $f (+\infty) = \lim_{x\to +\infty}f(x)\geq 0$. When the functio$x^r[f(x)-f(+\infty)]$ is completely monotonic on $(0,+\infty)$ if and only if $0 \leq r\leq \alpha$, the number $\alpha$, denoted by ${\rm deg}^x_{cm}[f]$, is called the completely monotonic degree of $f(x)$ with respect to $x\in (0,+\infty)$. The most integer $r$ which satisfies this condition is the completely monotonic integer degree of $f$. For more studies on complete monotonicity, the reader is also referred to \cite{guo1, kou2,kou3, kou6, Qi2, Qi3}.

For $x > 0$, the classical gamma function $\displaystyle\Gamma(x)=\int_0^\infty t^{x-1}e^{-t}dt$ first introduced by L. Euler, is one of the most important functions in mathematical analysis. It often appears in asymptotic series, hypergeometric series, Riemann zeta function, number theory, and so on.

In this section we give an answer to a problem suggested by Qi in his paper \cite{qi0}. Which states, for all $m\geq 0$, the function $(-1)^mx^m\Phi^{(m)}(x)$ is completely monotonic, where $\Phi(x)=x\psi'(x)-1$. We show that this is not true in general. Also, we investigate the $q$-analog of this problem. 

Before proving our results, let us recall the following facts.
For $m\geq 1$, we define on $(0,+\infty)$ the function
$$f_m(t)=\frac{d^{m-1}}{dt^{m}}\Big(\frac{t^{m}}{1-e^{-t}}\Big).$$
Recall the Hardy-Littlewood  entire function is defined on the complex plane as follows 
$\displaystyle H(z)=\sum_{k=1}^\infty\frac1k\sin(\frac z{k}),$ see \cite{alz1}. We set also $$s(z)=\displaystyle\frac12+\frac1\pi H(\frac z{2\pi}).$$ %By a simple argument of uniform convergence for series, we have for all $z\in\Bbb C$, and $p\geq 1$
%$$H^{(p)}(z)= \sum_{k=1}^\infty\frac1{k^{p+1}}\sin(\frac zk+\frac{p\pi}2).$$
In \cite{alz1} theorem 1.1, it is proved that,
$$\lim_{m\to\infty}\frac1{m!}f_m(z/m)=s(z),$$
the convergence begin uniform on every compact set of the complex plane. %Following the same method one proves
%\begin{proposition} For $z\in\Bbb C$, and all $p\in\Bbb N$
%$$\lim_{m\to\infty}\frac1{m!}f^{(p)}_m(z/m)=s^{(p)}(z),$$
%the convergence begin uniform on every compact set of the complex plane.
%\end{proposition}
Recall the result du to Alzet et al. See for instance \cite{alz1}.
\begin{proposition}\label{pr} There are positive constants $C$ and a positive sequence $x_{j_k}$ with $\lim_{k\to\infty}x_{j_k}=+\infty$ such that for all $k\in\Bbb N$,
$$H(x_{j_k})\leq -C\sqrt{\log\log x_{j_k}}.$$
\end{proposition}
\begin{proposition}\label{p01}
Let $$\Phi(x)=x\psi'(x)-1.$$ There exists $m_0\in\Bbb N$, such that for every $m\geq m_0$ and $\alpha\in\{m,m-1\}$, the functions $(-1)^mx^\alpha\Phi^{(m)}(x)$ are not completely monotonic.
\end{proposition}
{\bf Proof.} Using the Liebneitz rule we get for $m\geq 1$
$$\Phi^{(m)}(x)=x\psi^{(m+1)}(x)+m\psi^{(m)}(x),$$
Since, $$\psi^{(n)}(x)=(-1)^{n+1}\int_0^\infty\frac{t^n}{1-e^{-t}}e^{-x t}dt.$$
Thus,
$$(-1)^mx^m\Phi^{(m)}(x)=\int_0^\infty \frac{ t^{m+1}}{1-e^{-t}}x^{m+1}e^{-x t}dt-m\int_0^\infty \frac{t^{m}}{1-e^{-t}}x^me^{-x t}dt,$$
Integrate by part yields,
$$(-1)^mx^m\Phi^{(m)}(x)=\int_0^\infty (f_{m+1}(t)-mf_m(t))e^{-x t}dt,$$
where $f_m(t)=\Big(t^m/(1-e^{-t})\Big)^{(m)}$.

Since,
%$$f_{m+1}(t)=t f_m'(t)+(m+1)f_m(t).$$
%$$2(m+1)f_{m+1}-f_{m+2}-m(m+1)f_m$$
%$$=2(m+1)f_{m+1}-tf'_{m+1}-(m+2)f_{m+1}-m(m+1)f_m$$
%$$=mf_{m+1}-tf'_{m+1}-m(m+1)f_m$$
%$$=mtf'_{m}-tf'_{m+1}$$
%$$=mtf'_{m}-tf'_{m+1}$$
%$$=-t(2f_m'(t)+tf''_m)=-t(tf_m)''$$
Therefore, $$(-1)^mx^m\Phi^{(m)}(x)=\int_0^\infty g_m(t)e^{-x t}dt,$$
where $$g(t)=tf'_m(t)+f_m(t)=(tf_m(t))'.$$
If $(-1)^m\Phi^{(m)}(x)$ is completely monotonic for all $m$ and $t>0$. Then, $tf_m(t)$ increases on $(0,+\infty)$, and $f_m(t)\geq 0$. Moreover, By a result of Alzer \cite{alz1}, we saw that $\lim_{m\to\infty}(1/m!)f_m(t/m)=s(t)$ for all $t>0$. Which is impossible, due to the result Proposition \ref{pr}.

Furthermore,
$$(-1)^mx^{m-1}\Phi^{(m)}(x)=\int_0^\infty tf_m(t) e^{-x t}dt,$$
Using again Proposition \ref{pr}, we get the desired result. Namely, there is $m_0$, such that for $m\geq m_0$ the functions $(-1)^mx^{m-1}\Phi^{(m)}(x)$ and $(-1)^mx^{m}\Phi^{(m)}(x)$ are not completely monotonic.

%This gives, $(-1)^mx^m\Phi^{(m)}(x)$ is completely monotonic and then, ${\rm deg}_{cm}^x((-1)^m\Phi^{(m)}(x))\geq m$.
\begin{proposition} For all $m\in\Bbb N$, the function $(-1)^mx^{m+1}\Phi^{(m)}(x)$ is not completely monotonic.
\end{proposition}

{\bf Proof.} We have,
$$x\psi'(x)=\int_0^\infty h(t)e^{-x t}dt,$$
where, $h(t)=e^{-t}(e^t-1-t)/(1-e^{-t})^2=\Big(t/(1-e^{-t})\Big)'$
$$(-1)^mx^{m+1}\Phi^{(m)}(x)=\int_0^\infty h(t)t^mx^{m+1}e^{-x t}dt.$$
Therefore,
$$\Big((-1)^mx^{m+1}\Phi^{(m)}(x)\Big)^{(m+1)}=(m+1)!\int_0^\infty h(t)t^mL_{m+1}(xt)e^{-x t}dt.$$
Moreover,
$$t^mh(t)=t^m+\sum_{k=1}^\infty(1-k t)t^me^{-k t}.$$
Then,
$$\begin{aligned}\Big((-1)^mx^{m+1}\Phi^{(m)}(x)\Big)^{(m+1)}=&(m+1)!\int_0^\infty t^mL_{m+1}(xt)e^{-x t}dt\\&+(m+1)!\sum_{k=1}^\infty\int_0^\infty (t^m-kt^{m+1})L_{m+1}(xt)e^{-(x+k) t}dt\end{aligned}.$$
So, by the substitution $u=(x+k)t$, we get

$$\Big((-1)^mx^{m+1}\Phi^{(m)}(x)\Big)^{(m+1)}=(m+1)!\sum_{k=1}^\infty\int_0^\infty \Big(\frac{u^m}{(x+k)^{m+1}}-k\frac{u^{m+1}}{(x+k)^{m+2}}\Big)L_{m+1}(\frac x{x+k}u)e^{-u}du.$$
Moreover, $|L_{m+1}(\frac x{x+k}u)|\leq e^{x u/2(x+k) },$ so, for $x>0, m\geq 1$ then $|L_{m+1}(\frac x{x+k}u)| e^{-u}\leq e^{-u/2}$, and
$\displaystyle\lim_{x\to 0}\sum_{k=1}^\infty\frac1{(x+k)^m}L_{m+1}(\frac x{x+k}t)= \zeta(m)$. By the dominated convergence theorem, it follows that for $m\geq 1$,
$$\lim_{x\to 0^+}\Big((-1)^mx^{m+1}\Phi^{(m)}(x)\Big)^{(m+1)}=(m+1)!(m!\zeta(m+1)-(m+1)!\zeta(m+1))=-m(m+1)!m!\zeta(m+1).$$
So, $(-1)^{m+1}\Big((-1)^mx^{m+1}\Phi^{(m)}(x)\Big)^{(m+1)}$ is negative for $m$ odd. And the function $(-1)^mx^{m+1}\Phi^{(m)}(x)$ is not completely monotonic for $m$ odd.

On the other hand by differentiation we get
$$\begin{aligned}&\Big((-1)^mx^{m+1}\Phi^{(m)}(x)\Big)^{(m+2)}\\&=(m+1)!\sum_{k=1}^\infty\int_0^\infty \Big(-(m+1)\frac{u^m}{(x+k)^{m+2}}+k(m+2)\frac{u^{m+1}}{(x+k)^{m+3}}\Big)L_{m+1}(\frac x{x+k}u)e^{-u}du\\&+(m+1)!\sum_{k=1}^\infty\int_0^\infty \Big(\frac{u^m}{(x+k)^{m+1}}-k\frac{u^{m+1}}{(x+k)^{m+2}}\Big)(-\frac{ku}{(x+k)^2})L'_{m+1}(\frac x{x+k}u)e^{-u}du.\end{aligned}$$
Since, $L'_{m+1}(0)=-(m+1)$. It follows that
$$\begin{aligned}\lim_{x\to 0^+}\Big((-1)^mx^{m+1}\Phi^{(m)}(x)\Big)^{(m+2)}&=
(m+1)!\sum_{k=1}^\infty\int_0^\infty \Big(-(m+1)\frac{u^m}{k^{m+2}}+(m+2)\frac{u^{m+1}}{k^{m+2}}\Big)e^{-u}du\\&+(m+1)!\sum_{k=1}^\infty\int_0^\infty \Big(\frac{u^m}{k^{m+1}}-\frac{u^{m+1}}{k^{m+1}}\Big)(m+1)\frac{u}{k}e^{-u}du.\end{aligned}$$
$$\lim_{x\to 0^+}\Big((-1)^mx^{m+1}\Phi^{(m)}(x)\Big)^{(m+2)}=((m+1)!)^2(\zeta(m+2)-(m+2)\zeta(m+2)+\zeta(m+1)-(m+2)\zeta(m+1)).$$
$$\lim_{x\to 0^+}\Big((-1)^mx^{m+1}\Phi^{(m)}(x)\Big)^{(m+2)}=-((m+1)!)^2(m+1)(\zeta(m+2)+\zeta(m+1)).$$
So,
$$\lim_{x\to 0^+}(-1)^{m+2}\Big((-1)^mx^{m+1}\Phi^{(m)}(x)\Big)^{(m+2)}=(-1)^{m+1}((m+1)!)^2(\zeta(m+2)+\zeta(m+1)).$$
The right hand side is negative for $m$ even. Thus, the function $(-1)^mx^{m+1}\Phi^{(m)}(x)$ is not completely monotonic for $m$ even.
Which implies that the function $(-1)^mx^{m+1}\Phi^{(m)}(x)$ is not completely monotonic for all $m\geq 0$. Thus,
${\rm deg}^x_{cm}((-1)^m\Phi^{(m)}(x))<m+1$.
\begin{proposition} For all $m\in\Bbb N$, the function $(-1)^mx^{m-2}\Phi^{(m)}(x)$ is completely monotonic.
\end{proposition}
{\bf Proof.} Following the proof of Proposition \ref{p01}, we have
$$(-1)^mx^{m-2}\Phi^{(m)}(x)=\int_0^\infty \Theta_m(t) e^{-x t}dt,$$
where $$\Theta_m(t)=\int_0^tuf_m(u)du.$$
We saw that $$f_{m}(t)=m!+m!\sum_{k=1}^\infty e^{-kt}L_{m}(kt),$$
and by the fact that, $|L_m(x)|\leq e^{\frac x2}$ for all $x>0$. It then follows that
$$|\sum_{k=1}^\infty e^{-kt}L_{m-1}(kt)|\leq \frac{e^{-\frac t2}}{1-e^{-\frac t2}},$$
hence, $$m!(1-\frac{e^{-\frac t2}}{1-e^{-\frac t2}})\leq f_m(t),$$
It is easy seeing that $1-\frac{e^{-\frac t2}}{1-e^{-\frac t2}}\geq\frac12$ if and only if $t\geq 2\log 3$. Then, for all $ t\geq 2\log 3\simeq 2.19$
$$ f_m(t)\geq\frac{m!}2 t.$$
Moreover, H. Alzer et al. \cite{alz1} (p.113) showed that for all $u\in(-2\pi,2\pi)$,
$$f_m(u)=\int_0^\infty s(u x)x^{m}e^{-x}dx,$$
where, $$s(x)=\frac12+\frac1\pi\sum_{k=1}^\infty\frac1k\sin(\frac x{2k\pi}).$$
Let $t\in[0,2\pi)$, then,
$$\int_0^{t}uf_m(u)du=\int_0^{t}(\int_0^\infty s(u x)x^{m}e^{-x}dx)udu.$$
then,
$$\begin{aligned}\int_0^{t}uf_m(u)du&=\int_0^\infty\Big(\frac {t^2}4+\frac1\pi\sum_{k=1}^\infty\frac1k \int_0^{t}u\sin(\frac {x u}{2k\pi})du\Big))x^{m}e^{-x}dx\\&=2 \int_0^\infty\Big(\frac {t^2}8+\sum_{k=1}^\infty\left(\sin(\frac{t x}{2 k \pi })-\frac{t x}{2\pi k} \cos(\frac{t x}{2 k \pi }) \right)\Big)x^{m-2}e^{-x}dx.\end{aligned}$$
$$2 \int_0^\infty\Big(\sum_{k=1}^\infty\left(\frac {3(tx)^2}{(2 k\pi)^2}+\sin(\frac{t x}{2 k \pi })-\frac{t x}{2\pi k} \cos(\frac{t x}{2 k \pi }) \right)\Big)x^{m-2}e^{-x}dx.$$
Since, the function $3x^2-x\cos x+\sin x$ increases on $[0,+\infty)$ and then $3x^2-x\cos x+\sin x\geq 0$. Therefore, for all $t\geq 0$, $\displaystyle\int_0^{t}uf_m(u)du\geq 0$. This completes the proof.

\begin{proposition} For all $m\in\Bbb N$, the function $(-1)^mx^{m+1}\Phi^{(m)}(x)$ decreases on $(0,+\infty)$. The double inequality
$$\frac{m!}{2x^{m+1}}<(-1)^m\Phi^{(m)}(x)<\frac{m!}{x^{m+1}}.$$
is valid on $(0,+\infty)$ and sharp in the sense that the lower and upper bounds
cannot be replaced by any larger and smaller numbers respectively.
\end{proposition}
{\bf Proof.} Applying Liebnitz rule we get $m\geq 1$
\begin{equation}\label{e0}\Phi^{(m)}(x)=x\psi^{(m+1)}+m\psi^{(m)}(x),\end{equation}
Differentiate the function $x^{m+1}\Phi^{(m)}(x)$ yields
\begin{equation}\label{e}\Big((-1)^mx^{m+1}\Phi^{(m)}(x)\Big)'=x^m(-1)^m\Big(m(m+1)\psi^{(m)}+2(m+1)x\psi^{(m+1)}(x)+x^2\psi^{(m+2)}(x)\Big).\end{equation}
Using the integral representation
$$\psi^{(n)}(x)=(-1)^{n+1}\int_0^\infty\phi_n(t)e^{-x t}dt,$$
where $\phi_n(t)=t^n/(1-e^{-t})$.
It follows, that
$$\begin{aligned}&(-1)^m(m(m+1)\psi^{(m)}+2(m+1)x\psi^{(m+1)}(x)+x^2\psi^{(m+2)})(x)\\&=\int_0^\infty\Big(2(m+1)\phi_{m+1}'(t)
-\phi_{m+2}''(t)-m(m+1)\phi_m(t)\Big)e^{-x t}dt\end{aligned}.$$

By some algebra we get
$$2(m+1)\phi_{m+1}'(t)-\phi_{m+2}''(t)-m(m+1)\phi_m(t)=\frac{t^{1+m} e^t }{(e^t-1)^3} (e^t (2-t)-t-2).$$
Set $\Theta(t)=e^t (2-t)-t-2$, then $\Theta'(t)=e^t(1-t)-1$, and $\Theta''(t)=-t e^t$, moreover, $\Theta'(0)=0$ and $\Theta(0)=0$. Therefore, $\Theta(t)<0$ for all $t>0$. From equation \eqref{e}, one deduces that the function $(-1)^mx^{m+1}\Phi(x)$ is strictly increasing on $(0,+\infty)$. On the other hand, we saw that \begin{equation}\label{e1}\psi^{(m)}(x+1)=\psi^{(m)}(x)+\frac{(-1)^mm!}{x^{m+1}},\end{equation}
and for $x$ large enough
 \begin{equation}\label{e2}\psi^{(m)}(x)(x)=(-1)^{m-1}\Big(\frac{(m-1)!}{x^m}+\frac{m!}{2 x^{m+1}}+R_m(x)\Big),\end{equation}
with $\lim_{x\to\infty}x^{m+1}R_m(x)=0$.
%x\psi^{(m+1)}+m\psi^{(m)}(x)
By equation \eqref{e0}, we get
$$x^{m+1}\Phi^{(m)}(x)=x^{m+2}\psi^{(m+1)}(x+1)-(-1)^{m+1}(m+1)!+x^{m+1}\psi^{(m)}(x+1)-(-1)^{m}m(m)!,$$ and
$$\lim_{x\to 0}(-1)^mx^{m+1}\Phi^{(m)}(x)=m!.$$ Moreover, using \eqref{e2} we obtain
$$x^{m+1}\Phi^{(m)}(x)=(-1)^{m}(m!x +\frac{(m+1)!}{2 }+x^{m+2}R_{m+1}(x))+m(-1)^{m-1}((m-1)!x+\frac{m!}{2 x}+x^{m+1}R_m(x)).$$
Hence, $$\lim_{x\to +\infty}(-1)^mx^{m+1}\Phi^{(m)}(x)=\frac {m!}2.$$
This completes the proof.
\begin{proposition} For $q\in(0,1)$, let define on $(0,+\infty)$ the function
$$\Phi_q(x)=\frac{q^x-1}{\log q}\psi'_q(x)-q^x,$$
For all $n\geq 0$, the function $(-1)^n\Phi_q^{(n)}(x)$ is completely monotonic.

For $q>1$ the function $$\phi_q(x)=\frac{q^x-1}{q^x\log q}\psi'_q(x)-1,$$ is completely monotonic on $(0,+\infty)$.
\end{proposition}
{\bf Proof.} Since, $$\psi'_q(x)=\int_0^\infty\frac t{1-e^{-t}}e^{-xt}\gamma_q(dt),$$
where $\gamma_q$  is the positive discrete measure $$\gamma_q(t)=\left\{\begin{aligned}&-\log q\sum_{k=1}^\infty\delta_{t+k\log q}\;{\rm if}\;q\in(0,1)\\
&=t\;\;{\rm if}\;q=1\end{aligned}\right..$$
By differentiation with respect to $x$, we get
\begin{equation}\label{r1}\Big(\frac{q^x-1}{\log q} e^{-x t}\Big)^{(n)}=-\frac{1}{\log q}(-1)^nt^n e^{-x t}+\sum_{k=0}^nC_n^kq^x(\log q)^{k-1}(-1)^{n-k}t^{n-k} e^{-x t}.\end{equation}
Since, $$\Big(\frac{q^x-1}{\log q} e^{-x t}\Big)^{(n)}=-\frac{1}{\log q}(-1)^n(t^n-q^x(t-\log q)^n)e^{-x t} .$$
Thus, \begin{equation}\label{r}(-1)^n\Phi_q^{(n)}(x)=-\frac{1}{\log q}\int_0^\infty\frac t{1-e^{-t}}(t^n-q^x(t-\log q)^n)e^{-xt}\gamma_q(dt)-(-\log q)^n q^x.\end{equation}
%A direct induction gives, for all $m\geq 0$,
%(-1)^{n\Phi_q^{(n)}(x)=-\frac{1}{\log q}\int_0^\infty\frac t{1-e^{-t}}(t^n-q^x(t-\log q)^n)e^{-xt}\gamma_q(dt)-(-\log q)^n q^x
%One writes this integral as
%$$(-1)^n\Phi_q^{(n)}(x)=-\frac{1}{\log q}\int_0^\infty\frac t{1-e^{-t}}(t^n-q^x(t-\log q)^n)e^{-xt}\mu_q(dt),$$
%where $\mu_q$ is the positive measure $\mu_q=\gamma_q-(\log q)\delta_0$.

We split the integral on $(0,x_0)$ and $(x_0,+\infty)$ where $x_0=-(q^{\frac xn}\log q)/(1-q^{\frac xn})$.
Easy computations reveal that the function $U:t\mapsto t/(1-e^{-t})$ is strictly increasing on $(0,+\infty)$, with $\lim_{t\to 0} U(t)=1$. Let $n\geq 1$, and employing \eqref{r}.
%Let's $x_0=\frac{nq^x\log q}{q^x-1}$. Then, for $0<t<x_0$ we have, $U(t)< U(x_0)$ and if $t>x_0$, $U(t)> U(x_0)$ .
We get $$(-1)^n\Phi_q^{(n)}(x)\geq -\frac{U(x_0)}{\log q}\int_0^\infty (t^n -q^x(t-\log q)^n)e^{-x t}\mu_q(dt)-(-\log q)^n q^x.$$
Moreover, by using the fact that
$$\Theta_q(x):=\log q\frac{q^x}{q^x-1}=\int_0^\infty e^{-xt}\gamma_q(dt),$$
we get
$$ -\frac{1}{\log q}\int_0^\infty (t^n -q^x(t-\log q)^n)e^{-x t}\gamma_q(dt)=-\frac{(-1)^n}{\log q}\Theta_q^{(n)}(x)+\sum_{k=0}^nC_n^kq^x(\log q)^{k-1}(-1)^{k}\Theta_q^{(n-k)}(x).$$
$$ -\frac{1}{\log q}\int_0^\infty (t^n -q^x(t-\log q)^n)e^{-x t}\mu_q(dt)=-\frac{(-1)^n}{\log q}\Theta_q^{(n)}(x)+\frac{(-1)^n}{\log q}(q^x\Theta_q(x))^{(n)}=(-1)^n(\log q)^nq^x.$$
So, for all $n\geq 1$,
$$(-1)^n\Phi^{(n)}_q(x)\geq (-\log q)^nq^x(U(x_0)-1)\geq 0.$$
Furthermore, one deduces that $\Phi_q(x)$ decreases and we have $\lim_{x\to+\infty}\Phi_q(x)=0$ for $q\in(0,1)$. Thus, $\Phi_q(x)\geq 0$.

The second item follows from the fact that for $q>1$,
$\phi_q(x)=\Phi_{\frac1q}(x).$
 Which completes the proof.

\begin{corollary} The previous result is a generalization of the following one, the function $x\psi'(x)-1$ is completely monotonic. It corresponds to the case $q=1$.
\end{corollary}
In \cite{guo2}, it has been posed the following conjecture. The function $x^{x(\psi(x)-\log x)-\gamma},$ is logarithmically completely monotonic on $(0,+\infty)$, where $\gamma$ is the Euler constant. In the next proposition we prove a more general result.
%$$f(x)=x^{x(\psi(x)-\log x)-\frac\gamma{x \log x}}.$$
\begin{proposition}\label{p1}
The function $$f_\alpha(x)=x^{x(\psi(x)-\log x)-\alpha},$$ is logarithmically completely monotonic on $(0,+\infty)$ for $\alpha\geq-1/4$.
\end{proposition}
Recall the following results, see for instance \cite{sch}, [Theorem 3.6, p 19]
\begin{lemma}\label{lem1} Let $h:(0,+\infty)\rightarrow(0,+\infty)$ and $g$ be two functions.
\begin{enumerate}
\item If $h'$ and $g$ are completely monotonic, then $g\circ h$ is also completely monotonic.
\item If $(-\log(h))'$ is completely monotonic, then $h$ is also completely monotonic.
\end{enumerate}
\end{lemma}
{\bf Proof.} It is established in \cite{alz2} that the function $\theta_\beta(x)=x^\beta(\log x-\psi(x))$ is strictly completely monotonic on $(0,+\infty)$ if and only if $\beta\leq 1$. Moreover,
\begin{equation}\label{eq0}\theta_1(x)=\frac12+\int_0^\infty h(t)e^{-x t}dt,\end{equation}
where $h(t)=1/t^2-e^{-t}/(1-e^{-t})^2.$ In \cite{li}, it is proved that $h$ is strictly decreasing on $(0,+\infty)$ with $\lim_{x\to +\infty}h(t)=0$ and $\lim_{x\to 0}h(t)=1/12$. Hence, the function $\theta_1(x)-1/2$ is completely monotonic,
and, \begin{equation}\label{eq01}\frac12\leq\theta_1(x)\leq\frac12+\frac1{12x}.\end{equation}
Let $\varphi_\alpha(x)=(-\log f_\alpha(x))'$. Then,
\begin{equation}\label{eq1}\varphi_\alpha(x)=(\alpha\log x+\theta_1(x)\log x)',\end{equation}
and
\begin{equation}\label{eq2}\varphi_\alpha(x)=\frac{\alpha+\frac12}x+\frac1x(\theta_1(x)-\frac12)+\theta_1'(x)\log x\end{equation}
We show that $(-1)^n\varphi_\alpha^{(n)}(x)>0$ for $n\in\Bbb N$ and $x>0$.
%$$\varphi(x)=\frac\alpha x+\frac1x(\theta_1(x)-\frac12)+\theta_1'(x)\log x$$
Since, %$(-1)^n\Big(1/x(\theta_1(x)-\frac12)\Big)^{(n)}$
$$(-1)^{(n)}\varphi^{(n)}(x)=\frac{n!(\alpha+\frac12)}{x^{n+1}}+(-1)^n\theta_1^{(n+1)}(x)\log x +\sum_{k=1}^{n}C_n^k\frac{(k-1)!}{x^k}(-1)^{n-k+1}\theta_1^{(n-k+1)}(x)$$
$$+(-1)^n\Big(\frac1x(\theta_1(x)-\frac12)\Big)^{(n)}.$$
For every $x>0$
$$(-1)^n\Big(\frac1x(\theta_1(x)-\frac12)\Big)^{(n)}+\sum_{k=1}^{n}C_n^k\frac{(k-1)!}{x^k}(-1)^{n-k+1}\theta_1^{(n-k+1)}(x)\geq0,$$
and $(-1)^n\theta_1^{(n+1)}(x)\leq 0$. So, $(-1)^n\varphi^{(n)}_\alpha(x)\geq 0$ for $n\geq 0$ and $x\in(0,1)$

Let us assume, $x\geq 1$.
It remains to show that for $n\geq 1$, $x\geq 1$ and $\alpha\geq-1/4$, \begin{equation}\label{eq1}\frac{n!(\alpha+\frac12)}{x^{n+1}}+(-1)^n\theta_1^{(n+1)}(x)\log x\geq 0.\end{equation}
It sufficient to prove that for $n\geq 1$ and $x\geq 1$,
$$\frac{n!}{4x^{n+1}}+(-1)^n\theta_1^{(n+1)}(x)\log x\geq 0.$$

%Since, for $n\geq 1$, $(-1)^n\theta_0^{(n)}(x)=(-1)^{n+1}\psi^{(n)}(x)-(n-1)!/x^n$. Which gives the following double inequalities for the function $\theta_0$, for $n\geq 1$, $x>0$
%\begin{equation}\label{eq3}\frac{n!}{2x^{n+1}}\leq(-1)^n\theta_0^{(n)}(x)\leq\frac{n!}{x^{n+1}}.\end{equation}
%Using \eqref{eq1} and \eqref{eq3}. It sufficient to prove that for $n\geq 1$,
%$$\frac{n!}{2x^{n+1}}+(-1)^n\theta_1^{(n+1)}(x)\log x\geq 0.$$
On the first hand,
$$(-1)^n(\theta_1(x))^{(n+1)}\log x=(-1)^n(x\log x)^{(n+1)}\log x+(-1)^{n+1}(x\psi(x))^{(n+1)}\log x,$$
and $$(-1)^n(x\log x)^{(n+1)}=(-1)^n(\frac{(-1)^nn!}{x^n}+(n+1)\frac{(-1)^{n-1}(n-1)!}{x^n})=-\frac{(n-1)!\log x}{x^n}.$$
Then,
$$(-1)^n(\theta_1(x))^{(n+1)}\log x=((-1)^{n+1}(x\psi(x))^{(n+1)}-\frac{(n-1)!}{x^n})\log x,$$
Let us denote
$$g_n(x)=(-1)^n\theta_1^{(n+1)}(x)\log x+\frac{n!}{4x^{n+1}}.$$
$$g_n(x)=((-1)^{n+1}\Big(x\psi^{(n+1)}(x)+(n+1)\psi^{(n)}(x)\Big)-\frac{(n-1)!}{x^n})\log x+\frac{ n!}{4x^{n+1}},$$
$$\begin{aligned}&g_n(x+1)-g_n(x)\\
&=((-1)^{n+1}\Big((x+1)\psi^{(n+1)}(x+1)+(n+1)\psi^{(n)}(x+1)\Big)-\frac{(n-1)!}{(x+1)^n})\log (x+1)+\frac{ n!}{4(x+1)^{n+1}}\\
&-((-1)^{n+1}\Big(x\psi^{(n+1)}(x)+(n+1)\psi^{(n)}(x)\Big)-\frac{(n-1)!}{x^n})\log x-\frac{ n!}{4x^{n+1}}\\
&=((-1)^{n+1}\Big((x+1)\psi^{(n+1)}(x+1)+(n+1)\psi^{(n)}(x+1)\Big)-\frac{(n-1)!}{(x+1)^n})\log (1+1/x)\\
&+\Big((-1)^{n+1}\Big((x+1)\psi^{(n+1)}(x+1)+(n+1)\psi^{(n)}(x+1)-x\psi^{(n+1)}(x)
-(n+1)\psi^{(n)}(x)\Big)\\&-\frac{(n-1)!}{(x+1)^n}+\frac{(n-1)!}{x^n}\Big)\log x
+\frac14(\frac{ n!}{(x+1)^{n+1}}-\frac{ n!}{x^{n+1}})\\
&=((-1)^{n+1}\Big((x+1)\psi^{(n+1)}(x+1)+(n+1)\psi^{(n)}(x+1)\Big)-\frac{(n-1)!}{(x+1)^n})\log (1+1/x)\\
&+((-1)^{n+1}\Big(\frac{(-1)^{n+1}(n+1)!}{x^{n+1}}+\psi^{(n+1)}(x+1)
+\frac{(-1)^{n}(n+1)!}{x^{n+1}}\Big)-\frac{(n-1)!}{(x+1)^n}+\frac{(n-1)!}{x^n})\log x\\
&+\frac14(\frac{ n!}{(x+1)^{n+1}}-\frac{ n!}{x^{n+1}})\\
&=((-1)^{n+1}\Big((x+1)\psi^{(n+1)}(x+1)+(n+1)\psi^{(n)}(x+1)\Big)-\frac{(n-1)!}{(x+1)^n})\log (1+1/x)\\
&+\Big((-1)^{n+1}\psi^{(n+1)}(x+1)-\frac{(n-1)!}{(x+1)^n}+\frac{(n-1)!}{x^n}\Big)\log x+\frac14(\frac{ n!}{(x+1)^{n+1}}-\frac{ n!}{x^{n+1}}),
\end{aligned}$$
where we used in the third equation the following equation $$\psi^{(n+1)}(x+1)= \psi^{(n)}(x)\frac{(n-1)!}{x^{n}}+\frac{(-1)^nn!}{x^{n}},$$
Therefore, $$\begin{aligned}g_n(x+1)-g_n(x)&=((-1)^{n+1}((x+1)\psi(x+1))^{(n+1)}-\frac{(n-1)!}{(x+1)^n})\log (1+1/x)\\&+((-1)^{n+1}\psi^{(n+1)}(x+1)-\frac{(n-1)!}{(x+1)^n}+\frac{(n-1)!}{x^n})\log x\\
&+\frac14(\frac{ n!}{(x+1)^{n+1}}-\frac{ n!}{x^{n+1}})\end{aligned}$$ 
We will show that $g_n(x+1)-g_n(x)\leq 0$. Using the inequality due to Alzer for $x>0$ and $n\geq 1$
$$(-1)^{n+1}(x\psi(x))^{(n+1)}<\frac{(n-1)!}{x^n}.$$
It sufficient to prove that $$((-1)^{n+1}\psi^{(n+1)}(x+1)-\frac{(n-1)!}{(x+1)^n}+\frac{(n-1)!}{x^n})\log x+\frac14(\frac{ n!}{(x+1)^{n+1}}-\frac{ n!}{x^{n+1}})\leq 0.$$
%$$\leq\frac{(n-1)!}{x^n}\log(x+1)+((-1)^{n+1}\psi^{(n+1)}(x+1)-\frac{(n-1)!}{x^n})\log x$$
On the first hand,
$$(-1)^{n+1}\psi^{(n+1)}(x+1)-\frac{(n-1)!}{(x+1)^n}+\frac{(n-1)!}{x^n}=\int_0^\infty (1-e^{-t}-\frac{t^2}{e^t-1})t^{n-1}e^{-xt}dt.$$
Moreover, the derivative of the function $1-e^{-t}-\frac{t^2}{e^t-1}$ is given by $\frac{e^{-t}(1+te^t -e^t)^2}{(e^t-1)^2}$ and
$\displaystyle\lim_{t\to 0}1-e^{-t}-t^2/(e^t-1)=0$. Therefore, the function $(-1)^{n+1}\psi^{(n+1)}(x+1)-\frac{(n-1)!}{(x+1)^n}+\frac{(n-1)!}{x^n}$ is completely monotonic.

Let us denote, $$J(x)=\Big((-1)^{n+1}\psi^{(n+1)}(x+1)-\frac{(n-1)!}{(x+1)^n}+\frac{(n-1)!}{x^n}\Big)x
+\frac14(\frac{ n!}{(x+1)^{n+1}}-\frac{ n!}{x^{n+1}}),$$
and $$j(x)=((-1)^{n+1}\psi^{(n+1)}(x+1)-\frac{(n-1)!}{(x+1)^n}+\frac{(n-1)!}{x^n})\log x
+\frac14(\frac{ n!}{(x+1)^{n+1}}-\frac{n!}{x^{n+1}}),$$
By the previous remark, we get $j(x)\leq J(x)$ for $x>0$. On the other hand, and by using again the relation $n!/x^n=\int_0^\infty t^{n-1}e^{-xt}dt$, it follows that
$$J(x)=\int_0^\infty \Big(\frac{e^{-t}(1+te^t -e^t)^2}{(e^t-1)^2}+\frac14 t(e^{-t}-1)\Big)t^{n-1}e^{-xt}dt.$$
Let us denote $$\theta(t)=(1+te^t -e^t)^2-\frac14 t(e^t-1)^3.$$
By successive differentiation we get
$$\theta'(t)=\frac{1}{4} (1+e^{2 t}3-2 t+8 t^2+2 (-2+t) \cosh t+(2-8 t) \sinh t)),$$
$$\theta''(t)=\frac{1}{2} e^{2 t} (2+6 t+8 t^2-2 (1+t) \cosh t-(4+7 t) \sinh t)=\frac12e^{2t}\theta_1(t),$$
$$\theta_1'(t)=6+16 t-(6+7 t) \cosh t-(9+2 t) \sinh t,$$
Using the know inequalities $\cosh t\geq 1$ and $\sinh t\geq t$, we get $\theta_1'(t)\leq -2t^2<0$ for $t>0$. Moreover, $\theta_1(0)=0$
%$$\theta_2''(t)=-4+6 a e^t (11+3 t).$$
%One see that $\theta_2'$ is strictly positive for $t>0$ and $\theta_2(0)=3$.
Then, $\theta''$ is negative. Since, $\theta'(0)=0$ and $\theta(0)=0$, therefore, $\theta(t)\leq 0$, for $t\geq 0$% and $\theta'(0)=0$, $\theta(0)=0$. This implies that $\theta(t)< 0$ for $t>0$. Further,
$$J(x)=\int_0^\infty \frac{e^{-t}}{(e^t-1)^2}\theta(t)e^{-xt}dt.$$
Hence, $j(x)\leq J(x)< 0$. %Moreover, it is proved by Alzer, for $x>0$
%$$(-1)^{n+1}((x+1)\psi(x+1))^{(n+1)}-\frac{(n-1)!}{(x+1)^n}< 0.$$
Put all this together we get for all $x>0$, $g_n(x+1)< g_n(x)$ and by induction, for all $k\in\Bbb N$ and all $x>0$,
\begin{equation}\label{5}g_n(x+k)< g_n(x).\end{equation} Recall that
$$g_n(x)=((-1)^{n+1}\Big(x\psi^{(n+1)}(x)+(n+1)\psi^{(n)}(x)\Big)-\frac{(n-1)!}{x^n})\log x+\frac{ n!}{4x^{n+1}},$$
Using the double inequality  $$\frac{(n-1)!}{x^{n}}+\frac{n!}{2x^{n+1}}\leq(-1)^{n+1}\psi^{(n)}(x)\leq \frac{(n-1)!}{x^{n}}+\frac{n!}{x^{n+1}},$$ we get for all $n\geq 1$, $\lim_{x\to\infty}g_n(x)=0$. Letting $k$ to infinity in equation \eqref{5} we get $g_n(x)\geq 0$ for all $x\geq 1$. This proves
$(-1)^n\varphi^{(n)}(x)\geq 0$ for all $x\geq 1$ and all $n\geq 1$.
Since,
$$\varphi(x)=\frac\alpha x+\frac1x(\theta_1(x)-\frac12)+\theta_1'(x)\log x.$$
Moreover,
$$\theta_1(x)=\int_0^\infty h(t)e^{-xt}dt.$$
Differentiate yields
$$\theta'_1(x)=-\int_0^\infty t h(t)e^{-xt}dt=-\frac1x\int_0^\infty(th(t))'e^{-x t}dt.$$
Since, $(th(t))'=-1/t^2+1/(2-2 \cosh t)+(1/8) t\sinh t/\sinh^4(t/2)$ and one shows that $(th(t))'$ is bounded on $(0,+\infty)$ by a some positive constant $C$, say, indeed, $\lim_{t\to 0}(th(t))'=1/12$ and $\lim_{t\to\infty}(th(t))'=0$. Let $x>0$, then $|\theta'_1(x)|\leq C/x^2$ and $\lim_{x\to \infty}\theta'_1(x)\log x=0$.
Moreover, $\lim_{x\to\infty}\frac\alpha x+\frac1x(\theta_1(x)-\frac12)=0$ and $\varphi$ decreases on $(0,+\infty)$. All this together imply that $\varphi(x)\geq 0$ for all $x>0$. This completes the proof.
 
\begin{corollary} The function $$f_\alpha(x)=x^{x(\psi(x)-\log x)-\alpha},$$ is completely monotonic on $(0,+\infty)$ for $\alpha\geq-1/4$.
\end{corollary}
{\bf Proof.} Firstly, $f_\alpha(x)>0$. Applying Lemma \ref{lem1} and Proposition \ref{p1}, we get that  for $\alpha\geq -1/4$, $f_\alpha$ is completely monotonic.
\begin{corollary} For all $n\geq 2$, $x>1$, and $\alpha\geq-1/4$,
$$0<(-1)^{n}\theta_1^{(n)}(x)\leq \frac{(n-1)!(\alpha+\frac12)}{x^{n}\log x}$$
\end{corollary}
\begin{remark} It remains an open problem to prove that the function $f_\alpha(x)$ is completely monotonic on the range $(-1/4,-1/2]$
\end{remark}
%For the converse. If $f_\alpha$ is strictly completely monotonic, then $f_\alpha(x)>0$ and $f_\alpha'(x)< 0$. Whence, $(-\log f_\alpha(x))'< 0$, using the same proof as above we get $\alpha\geq-1/2$.


\begin{thebibliography}{9}
\bibitem{abr} M. Abramowitz and I. A. Stegun (Eds), Handbook of Mathematical Functions with Formu-
las, Graphs, and Mathematical Tables, National Bureau of Standards, Applied Mathematics
Series 55, 10th printing, Dover Publications, New York and Washington, 1972.
\bibitem{zeg} Abramowitz, M., Stegun, I. A., Handbook of mathematical functions, Dover, New York (1964).
\bibitem{alz2}  H. Alzer, On some inequalities for the gamma and psi functions, Math. Comp. 66 (1997),
no. 217, 373-389; available online at https://doi.org/10.1090/S0025-5718-97-00807-7.
\bibitem{alz1}  Alzer, H.; Berg, C.; Koumandos, S. On a conjecture of Clark and Ismail. J. Approx. Theory 2005, 134, 102-113.
  \bibitem{bi} Bierens de Haan, D., Nouvelles tables d'int\'egrales d\'efinies, Amsterdam, 1867. (Reprint) G. E.
Stechert \& Co., New York, 1939.

\bibitem{clar} Clark, W.E.; Ismail, M.E.H. Inequalities involving gamma and psi functions. Anal. Appl. 2003, 1, 129-140.
\bibitem{chen} C.-P. Chen, F. Qi, and H. M. Srivastava, Some properties of functions related to the gamma
and psi functions, Integral Transforms Spec. Funct. 21 (2010), no. 2, 153--164; available online
at https://doi.org/10.1080/10652460903064216.
\bibitem{et} Erd\'elyi, A. et al., Tables of Integral Transforms, vols. I and II. McGraw Hill, New York, 1954.
\bibitem{guo1}     Guo, B.-N., Qi, F.: A completely monotonic function involving the tri-gamma function and with degree one. Appl. Math. Comput. 218, 9890-9897 (2012). https://doi.org/10.1016/j.amc.2012.03.075.
 \bibitem{guo2}   B. N. Guo, Y. J. Zhang, F. Qi, Refinements and sharpenings of some double inequalities for
bounding the gamma function, J. Ineq. Pure Appl. Math. 9(1) (2008) Article 17.
    \bibitem{kou1} Koumandos, S. Remarks on some completely monotonic functions. J. Math. Anal. Appl. 2006, 324, 1458-1461.
 \bibitem{kou2} Koumandos, S.: Monotonicity of some functions involving the gamma and psi functions. Math. Compet. 77, 2261-2275 (2008). https://doi.org/10.1090/s0025-5718-08-02140-6
 \bibitem{kou3} Koumandos, S., Lamprecht, M.: Some completely monotonic functions of positive order. Math. Compet. 79, 1697-1707 (2010). https://doi.org/10.1090/s0025-5718-09-02313-8
 \bibitem{kou4} Koumandos, S., Lamprecht, M.: Complete monotonicity and related properties of some special functions. Math. Compet. 82, 1097-1120 (2013). https://doi.org/10.1090/s0025-5718-2012-02629-9
  \bibitem{kou5}Koumandos, S., Pedersen, H.L.: Completely monotonic functions of positive order and asymptotic expansions of the logarithm of Barnes double gamma function and Euler's gamma function. J. Math. Anal. Appl. 355, 33-40 (2009). https://doi.org/10.1016/j.jmaa.2009.01.042
  \bibitem{kou6} Koumandos, S., Pedersen, H.L.: Absolutely monotonic functions related to Euler's gamma function and Barnes' double and triple gamma function. Monatshefte Math. 163, 51-69 (2011). https://doi.org/10.1007/s00605-010-0197-9.
   \bibitem{kou7}    S. Koumandos, Remarks on some completely monotonic functions, J. Math. Anal. Appl. 324
(2006), no. 2, 1458-1461; available online at http://dx.doi.org/10.1016/j.jmaa.2005.12.
017.
\bibitem{li}  A.-Q. Liu, G.-F. Li, B.-N. Guo, and F. Qi, Monotonicity and logarithmic concavity of two
functions involving exponential function, Internat. J. Math. Ed. Sci. Tech. 39 (2008), no. 5,
686--691; Available online at http://dx.doi.org/10.1080/00207390801986841. 2
   \bibitem{olv}  F. W. J. Olver, D. W. Lozier, R. F. Boisvert, and C. W. Clark (eds.), NIST Handbook of
Mathematical Functions, Cambridge University Press, New York, 2010; available online at
http://dlmf.nist.gov/.
   \bibitem{Qi1}   Qi, F., Guo, B.-N.: Lévy-Khintchine representation of Toader-Qi mean. Math. Inequal. Appl. 21, 421-431 (2018).
https://doi.org/10.7153/mia-2018-21-29
\bibitem{qi0} Qi, F. Complete monotonicity and monotonicity of two functions defined by two derivatives of a function involving trigamma function.
 \bibitem{Qi2} Qi, F., Guo, B.-N.: The reciprocal of the weighted geometric mean of many positive numbers is a Stieltjes function. Quaest. Math. 41, 653-664 (2018). https://doi.org/10.13140/RG.2.2.23822.36163
 \bibitem{Qi3} Qi, F., Li, W.-H.: Integral representations and properties of some functions involving the logarithmic function. Filomat 30, 1659-1674 (2016). https://doi.org/10.2298/FIL1607659Q
 \bibitem{Qi4} Qi, F., Lim, D.: Integral representations of bivariate complex geometric mean and their applications. J. Comput. Appl. Math. 330, 41-58 (2018). https://doi.org/10.1016/j.cam.2017.11.047
   \bibitem{Qi5}   Qi, F., Liu, A.-Q.: Completely monotonic degrees for a difference between the logarithmic and psi functions. J. Comput. Appl. Math. 361, 366-371 (2019). https://doi.org/10.1016/j.cam.2019.05.001.
    \bibitem{Qi6}   F. Qi, A double inequality for the ratio of two non-zero neighbouring Bernoulli numbers, J.
Comput. Appl. Math. 351 (2019), 1-5; available online at https://doi.org/10.1016/j.cam.
2018.10.049.
 \bibitem{Qi7} F. Qi, Completely monotonic degree of a function involving the tri- and tetra-gamma functions, arXiv preprint (2013), available online at http://arxiv.org/abs/1301.0154
 \bibitem{Qi8} F. Qi, Notes on a double inequality for ratios of any two neighbouring non-zero Bernoulli
numbers, Turkish J. Anal. Number Theory 6 (2018), no. 5, 129-131; available online at
https://doi.org/10.12691/tjant-6-5-1.
 \bibitem{Qi9}  F. Qi and R. J. Chapman, Two closed forms for the Bernoulli polynomials, J. Number Theory
159 (2016), 89-100; available online at https://doi.org/10.1016/j.jnt.2015.07.021.

 \bibitem{Qi10} F. Qi, G. Bai-Ni. Completely monotonicities of functions involving the Gamma and Digamma functions, RGMIA Res. Rep. Coll. 7, no. 1 (2004) Art. 6.
      \bibitem{Qi11} F. Qi, G. Bai-Ni. Some completely monotonic functions involving the Gamma and Polygamma functions, RGMIA Res. Rep. Coll. 7, no. 1 (2004) Art. 8.
 \bibitem{mit} Mitrinov\'ic, D.S., Pe\'cari\'c, J.E., Fink, A.M.: Classical and New Inequalities in Analysis. Kluwer Academic, Dordrecht (1993).

\bibitem{sch}Schilling, R.L., Song, R., Vondra¡cek, Z.: Bernstein Functions-Theory and Applications, 2nd edn. de Gruyter Studies in Mathematics, vol. 37. de Gruyter, Berlin (2012).

\bibitem{wid}  Widder, D.V.: The Laplace Transform. Princeton University Press, Princeton (1946).

\bibitem{xu} Y. Xu and X. Han, Complete monotonicity properties for the gamma function and Barnes
G-function, Sci. Magna 5 (2009), no. 4, 47-51.
\bibitem{zeg} Abramowitz, M., Stegun, I. A., Handbook of mathematical functions, Dover, New York (1964).
\end{thebibliography}
\end{document}